\documentclass[onefignum,onetabnum]{siamart171218}

\usepackage{amsfonts}
\usepackage{xfrac}
\usepackage{graphicx}
\usepackage{epstopdf}
\usepackage{algorithmic}
\usepackage{amssymb}
\usepackage{booktabs}
\usepackage{algorithmic}
\Crefname{ALC@unique}{Line}{Lines}
\usepackage[T1]{fontenc}
\usepackage{tikz,caption}
\usepackage[caption=false,font=footnotesize]{subfig}
\usepackage{pgfplots}
\pgfplotsset{compat=newest}
\pgfplotsset{plot coordinates/math parser=false}
\usepackage{url}
\usepackage{color}
\usetikzlibrary{plotmarks} 
\usetikzlibrary{external}
\usepgfplotslibrary{external} 
\newcommand{\R}{\mathbb{R}}

\newcommand{\norm}[1]{\lVert#1\rVert}

\newcommand{\W}{{\ensuremath{{W}}}}

\newcommand{\V}{{V}}

\newcommand{\I}{{I}}
\newcommand{\Ll}{{L}}

\def\matlab{{\large {\sc matlab}}}

\def\norm#1{\left\|#1\right\|} 
\def\vec#1{\mathrm{vec}(#1)} 

\newtheorem{example}[theorem]{Example}

\definecolor{DarkBlue}{rgb}{0,0.08,0.45}
\definecolor{DarkRed}{rgb}{.65,0,0}

\newcommand{\TheTitle}{Solving differential Riccati equations: A nonlinear space-time method using tensor trains} 
\newcommand{\ShortTitle}{Solving DREs: A nonlinear space-time method} 
\newcommand{\TheAuthors}{Tobias Breiten and Sergey Dolgov and Martin Stoll}

\headers{\ShortTitle}{\TheAuthors}
\title{{\TheTitle}\thanks{}}

\author{Tobias Breiten\thanks{Institute of Mathematics and Scientific Computing, University of Graz, 8010 Graz, Austria ({\tt tobias.breiten@uni-graz.at })} \and Sergey Dolgov\thanks{University of Bath, Department of Mathematical Sciences, BA2 7AY Bath, United Kingdom ({\tt s.dolgov@bath.ac.uk})} \and Martin Stoll\thanks{Technische Universit\"at Chemnitz, Department of Mathematics, Scientific Computing Group, 09107 Chemnitz, Germany, ({\tt martin.stoll@mathematik.tu-chemnitz.de})}}

\begin{document}
\maketitle
\begin{abstract}
Differential algebraic Riccati equations are at the heart of many applications in control theory. They are time-depent, matrix-valued, and in particular nonlinear equations that require special methods for their solution. Low-rank methods have been used heavily computing a low-rank solution at every step of a time-discretization. We propose the use of an all-at-once space-time solution leading to a large nonlinear space-time problem for which we propose the use of a Newton--Kleinman iteration. Approximating the space-time problem  in low-rank form requires fewer applicatons of the discretized differential operator and gives a low-rank approximation to the overall solution.
\end{abstract}

\begin{keywords} 
Optimal Control, Low-rank methods, Riccati equations, Non-linear problems
\end{keywords}

\begin{AMS}
	65F15, 
	65T50, 
	68T05, 
	62H30  
\end{AMS}

\section{Motivation and main challenge} Linear and nonlinear matrix equations play an important role in control theory, see, e.g., \cite{AboFIJ03}. Well-known examples are the \emph{Lyapunov} and the \emph{Sylvester} equation which are intimately  related to stability, controllability, and observability concepts for linear dynamical control systems of the form 
\begin{equation}\label{eq:lti_sys}
\begin{aligned}
\dot{x}(t)&=Ax(t) + Bu(t), \ \ x(0)= x_0 , \\
y(t)&= Cx(t),
\end{aligned}
\end{equation}
where $A \in \mathbb R^{n\times n} , B\in \mathbb R^{n\times m}$ and $C\in \mathbb R^{p\times n}$. In the context of finite-horizon optimal feedback control of \eqref{eq:lti_sys}, one has to solve \emph{differential Riccati equations (DRE)}  
\begin{align}\label{eq:dre}
\dot{P} + A^T P +  PA  -PBB^TP + C^T C&= 0, \ \ P(t_f) = M,
\end{align}
where $M\in \mathbb R^{n\times n}$ is associated with a terminal penalty term at the final time horizon $t_f$, see the details below. From a computational point of view, equation \eqref{eq:dre} poses many numerical challenges. Indeed, in this case we have to solve for the unknown $P(\cdot)\colon [0,t_f]\to \mathbb R^{n\times n}$, i.e., a time-varying matrix with $n^2$  entries in each step. If the system \eqref{eq:lti_sys} results from a spatial semi-discretization of a partial differential equation (PDE), a computationally efficient method is crucial. A key ingredient in this context is that of a \emph{low numerical rank}. For the algebraic counterpart of \eqref{eq:dre}, it is known, see, e.g., \cite{BenS13,LinS15,Opm15,Sim16} that $P\approx LL^T$ can be well approximated by low-rank factors $L\in \mathbb R^{n\times k}$ with $k\ll n$ if at least the control or the observation matrices $B$ and $C$ correspond to finite-dimensional operators. In particular, this holds true for many relevant PDEs of parabolic type (cf. \cite{Opm15}) and practically realisable controllers.

In the time-varying case, only a few theoretical results have been obtained. For a recent discussion on low-rank solutions of differential Riccati equations, we refer to \cite{Sti18}. Nevertheless, many numerical approaches exist and most of them rely on a time discretization of \eqref{eq:dre} for an abstract nonlinear matrix-valued ordinary differential equation (ODE) of the form
\begin{align*}
\dot{P}(t) = f(t,P(t)) , \ \ P(0)=P_0.
\end{align*}
Here, we consider \eqref{eq:dre} as running forward in time. In order to avoid unnecessarily small step sizes due to stiff problems, many authors have proposed to use implicit schemes such as, e.g., Peer methods \cite{LanMS15,Lan18}, Rosenbrock methods \cite{BenM13} or BDF methods \cite{BenM04}. A common idea for all these methods now is that on a given time grid $0=t_0,t_1,\dots,t_{n_t}=t_f$ a discrete approximation $\widehat{P}$ of $P$ is obtained by solving nonlinear equations of the form
\begin{align}\label{eq:abstract_implicit}
P_{i+1} = G(P_{i+1-s},\dots,P_{i}) + H(t_{i+1},P_{i+1})
\end{align}
where the operators $G$ and $H$ depend on the chosen time stepping scheme with $s$ stages.
In the particular case of the DRE, an implicit Euler scheme for example results in a series of algebraic Riccati equations
\begin{align*}
\left(\tau A-\frac{1}{2}I\right)^\top P_{i} + P_{i} \left(\tau A-\frac{1}{2}I\right) - \tau P_i BB^\top P_i + \tau C^\top C + P_{i+1}= 0,  
\end{align*}
with $\tau$ denoting the time step size. Note that this equation runs backward in time with $P_{i+1}$ already known at step $i$.
For such equations, efficient numerical methods that are based on a low-rank Newton-Kleinman iteration (more details are given below) \cite{MenOPP18,FeiHS09} are known to perform well. 

A slightly different approach is given by splitting methods \cite{HanS14} that divide \eqref{eq:dre} into two parts 
\begin{align*}
\dot{P}(t) &= F(P(t)) + G(P(t)), \ \ P(0)=P_0 \\
F(P(t))&= A^TP(t) + P(t)A+C^TC , \\
G(P(t))&= -P(t)BB^TP(t),
\end{align*}
and, instead of \eqref{eq:dre}, rather consider the two individual problems 
\begin{align*}
\dot{P}(t) &= F(P(t)), \ \ P(0)=P_0, \\
\dot{P}(t) &= G(P(t)), \ \ P(0)=P_0.
\end{align*}
The benefit of this decomposition is that the first equation is affine and the second equation can be solved exactly, see \cite[Lemma 3.2]{HanS14}. One can thus apply splitting schemes such as, e.g., Lie splitting or Strang splitting. As for the time stepping schemes, the key idea is to approximate the unknown $P(\cdot)$ for each time step $t_i$ via low-rank factors. 

Let us emphasize that the previous methods all rely on storing the information for each time step separately. Hence, if each matrix $P(t_i)\approx L_iL_i^T$, with $L_i \in \mathbb R^{n\times k_i}$ one still has to store $n\cdot \sum_{i=1}^{t_f}k$ entries. While this is reasonable for problems with few time steps, it might cause difficulties if long time horizons $t_f$ or fine time grids have to be considered. The strategy we pursue in this manuscript is different in the sense that we proceed by discretizing the DRE in space and time simultaneously.  Of course, we then have a nonlinear problem of vast dimensionality defined on the space-time cylinder.
We here propose the use of multivariate tensor decompositions that allow us to easily approximate the solution of a linear problem defined over a high-dimensional space. This means that we propose the use of an outer nonlinear solver of Newton-type that at its core needs to solve a linearized problem. For this we propose low-rank tensor methods that depend on the rank $r$ of the solution over the whole time-interval and assumes this rank to be (approximately) constant. As our numerical examples confirm, this seems to be a reasonable assumption and could lead to a significant reduction in assembly time and storage requirements. \\

\textbf{Structure of the paper.}
In the next section we introduce the main idea of low-rank methods for all-at-once approaches. These methods tackle the solution of the matrix equation in a holistic way, where both space and time are discretized simultaneously and then solved in a coupled fashion. We also motivate that the rank of a space-time-solution is often found to be small so that our approach seems a feasible alternative. The heart of this paper is given in Section~\ref{sec::lrdre}, where we introduce the low-rank method for solving the differential Riccati equation. The method of choice relies on the tensor train format, which we introduce in Section~\ref{sec::tt}. In the last section we present results for various setups and illustrate that our suggested all-at-once scheme provides a viable alternative to other methods.

\textbf{Notation.}
We start by recalling several common notation from numerical\- (multi)\-linear algebra. Given two matrices $V \in \mathbb R^{p \times q}$ and $W \in \mathbb R^{n\times m}$, the \emph{Kronecker product} is defined by 

$$
W\otimes V:=
\left[
\begin{array}{ccc}
w_{11}V&\hdots&w_{1m}V\\
\vdots&\ddots&\vdots\\
w_{n1}V&\hdots&w_{nm}V\\
\end{array}
\right].
$$
We further introduce the \emph{vectorization operator}
\begin{align*}
\mathrm{vec}&\colon \mathbb R^{n\times m} \to \mathbb R^{nm}, \quad
\vec{W} = \vec{[w_1,\dots,w_m]} := \begin{bmatrix} w_1\\ \vdots \\ w_m \end{bmatrix},
\end{align*}
and remind the reader of the following relation between these two operators
\begin{align*}
\left(W^T\otimes V\right)\vec{Y}=\vec{VYW}.
\end{align*}


\section{Low-rank methods for all-at-once problems: motivation via open--loop control}
\label{sec::lrol}
The \textit{curse of dimensionality} \cite{Bel13} is encountered by researchers from many different communities such as large-scale PDE-constrained optimization \cite{StoB15}, uncertainty quantification \cite{SchS13,Zhaetal15}, chemical engineering \cite{FuP18}, statistics \cite{Donetal00} and others. Facing the optimal control problem given above, we rely on a technique introduced by the authors in \cite{StoB15}. Their main idea is to discretize the optimization problem in space and time in an all-at-once fashion. This of course leads to a very large-dimensional system. In \cite{StoB15}, the focus was on computing an optimal open-loop control for a linear system via solving the associated optimality system. In that case, one is faced with a large linear system of equations which, due to the highly structured system matrix, can be efficiently handled by iterative low-rank solvers. For several PDEs, this allows to reduce the necessary storage requirements to a fraction of the original problem dimension. Such ideas have recently been introduced for several applications and are in part due to the introduction of highly sophisticated compressed tensor formats \cite{GraKT13,Ose11} and correspondingly elegant solvers \cite{Dol13,DolS14,GraKT13,HolRS12}.   

Since our approach for \eqref{eq:dre} crucially depends on the above ideas, we provide a more detailed introduction into the topic on the example of open-loop control. For this purpose, let us consider 
a control system of the form
\begin{align*}
E\dot{x}(t) &= A x(t) + Bu(t), \ \ x(0) = x_0 , \\
y(t)&= Cx(t), 
\end{align*}
where $E,A \in \mathbb R^{n\times n}$,  $B\in \mathbb R^{n\times m}$ denotes a control operator and $C\in \mathbb R^{p\times n}$ is an observation operator. Given a \emph{desired trajectory} $x_d\colon [0,t_f] \to \mathbb R^n$, we want to solve the optimal control problem 
\begin{equation}
\label{jyu2}
\begin{aligned}
\min_{\substack{x\in L^2(0,t_f;\mathbb R^n)\\u\in L^2(0,t_f;\mathbb R^m)}}
J(x,u):&=\frac{1}{2}\int\limits_0^{t_f} \| y(t) - Cx_d(t) \|^2 \, \mathrm{d}t+\frac{\beta}{2}\int\limits_0^{t_f}\| u(t)\|^2 \,  \mathrm{d}t + \frac{1}{2} \| y(t_f)- Cx_d(t_f) \|^2 \\
\text{s.t. } \  E\dot{x}(t)&= Ax(t) + Bu(t), \ \ x(0)=x_0.
\end{aligned}
\end{equation}
Well-known results from optimal control theory, see, e.g., \cite{Tro05} now imply that the optimal solution $(\bar{x},\bar{u})$ is determined by an optimality system of the form
\begin{equation}\label{eq:cont_opt_sys}
\begin{aligned}
E\dot{\bar{x}}(t)&= A\bar{x}(t) + B\bar{u}(t) , && \bar{x}(0)=x_0, \\
-E^T\dot{p}(t) &= A^T p(t) + C^TC (\bar{x}(t)-x_d(t)), && p(t_f) = C^TC ( \bar{x}(t_f)- x_d(t_f)) , \\
0&=\beta \bar{u}(t)+ B^T p(t)  , && 
\end{aligned}
\end{equation}
where $p$ denotes the adjoint state. The approach proposed in \cite{StoB15} now considers the (time) discrete version of the previous optimality system. For this purpose, assume that the discretized state, adjoint, and control are described by the vectors 
\begin{align*}
\mathbf{x} = \begin{bmatrix} x_1 \\ \vdots \\ x_{n_t} \end{bmatrix}\in \mathbb R^{n\times n_t}, \quad 
\mathbf{u} = \begin{bmatrix} u_1 \\ \vdots \\ u_{n_t} \end{bmatrix}\in \mathbb R^{m\times n_t}, \quad
\mathbf{p} = \begin{bmatrix} p_1 \\ \vdots \\ p_{n_t} \end{bmatrix}\in \mathbb R^{n\times n_t}.
\end{align*}
Using the rectangle rule for the cost functional $J$ and the implicit Euler scheme for the control system, we obtain the following discrete analogue of \eqref{eq:cont_opt_sys}
\begin{equation}\label{eq:disc_opt_sys}
\begin{aligned} 
\underbrace{
	\left[
	\begin{array}{ccc}
	I_{n_t}\otimes\tau C^TC&0&-\left(\I_{n_t}\otimes L^T+G^{T}\otimes E^T\right)\\
	0&I_{n_t}\otimes\beta\tau I_m &I_{n_t}\otimes\tau B^T\\
	-\left(\I_{n_t}\otimes\Ll+G\otimes E\right)&I_{n_t}\otimes\tau B&0\\
	\end{array}
	\right]
}_{\mathcal{A}}
\left[
\begin{array}{c}
\mathbf{x}\\
\mathbf{u}\\
\mathbf{p} \\
\end{array}
\right]=\mathbf{f},
\end{aligned}
\end{equation}
where $L=E-\tau A$ and $G$, due to the implicit Euler scheme, is of the form 
\begin{align*}
G = \begin{bmatrix} 0 & & & \\ -1 & 0 & & \\ & \ddots & \ddots & \\ & & -1 & 0 \end{bmatrix} \in \mathbb R^{n_t \times n_t}
\end{align*} 
and the right hand side $\mathbf{f}$ is determined by $x_0$ and $x_d$, respectively. 

With the intention of reducing the storage needed for the vectors $\mathbf{x},$ $\mathbf{u},$ and $\mathbf{p}$, let us rather consider matrix representations of the form 
\begin{align*}
X=[x_1,\dots,x_{n_t}] \in \mathbb R^{n\times n_t} , \ \ 
U=[u_1,\dots,u_{n_t}] \in \mathbb R^{m\times n_t} , \ \ 
P=[p_1,\dots,p_{n_t}] \in \mathbb R^{n\times n_t},
\end{align*} 
which we approximate by low-rank representations
\begin{equation}\label{eq:low_rank_decomp}
\begin{aligned}
X&\approx\W_{{X}}\V_{{X}}^{T}\textnormal{ with }
\W_{{X}}\in\R^{n\times k_1},\V_{{X}}\in\R^{n_t\times k_1}\\
U&\approx\W_{{U}}\V_{{U}}^{T}\textnormal{ with }
\W_{{U}}\in\R^{m\times k_2},\V_{{U}}\in\R^{n_t\times k_2}\\
P&\approx\W_{{P}}\V_{{P}}^{T}\textnormal{ with }
\W_{{P}}\in\R^{n\times k_3},\V_{{P}}\in\R^{n_t\times k_3}
\end{aligned}
\end{equation}
with $k_{1,2,3}$ being small in comparison to $n_t.$ One can then implement the iterative solver to maintain the low-rank style of the solution in combination with an additional truncation scheme based on a truncated SVD or a QR reduction \cite{StoB15,KreT10,BenOS13}.  

\newlength\fheight
\newlength\fwidth
\setlength\fheight{0.25\linewidth} 
\setlength\fwidth{0.31\linewidth} 

\begin{example}
	Let us consider the one-dimensional heat equation 
	\begin{align*} 
	\frac{\partial}{\partial t} x(\xi,t) &= \frac{\partial^2}{\partial \xi^2} x(\xi,t) + \chi_{\omega} u(\xi,t) && \text{in } (0,1) \times (0,2), \\
	\frac{\partial}{\partial \xi} x(0,t)&=0=\frac{\partial}{\partial \xi} x(1,t) && \text{in } (0,2), \\
	x(\xi,0)&=\frac{1}{\sqrt{ 0.05 \pi}} \mathrm{exp}\left(-\frac{(\xi-0.25)^2}{0.05}\right) && \text{in } (0,1), 
	\end{align*}
	with control domain $\omega=(0.1,0.4)\cup (0.6,0.9)$ and desired trajectory (see Fig.~\ref{fig:motivation_1d_heat}, middle) 
	\begin{align*}
	x_d(\xi,t)= (t-2)x(\xi,0)+t\frac{1}{\sqrt{ 0.02 \pi}} \mathrm{exp}\left(-\frac{(\xi-0.75)^2}{0.02}\right).
	\end{align*}
	In Figure~\ref{fig:motivation_1d_heat} (left), we show the actual state trajectory obtained by solving the discrete optimality system \eqref{eq:disc_opt_sys} corresponding to a spatio-temporal discretization with $n=1000$ and $n_t=2000$ grid points.
	\begin{figure}[htb]
		\begin{center}
			\includegraphics[width=0.49\linewidth]{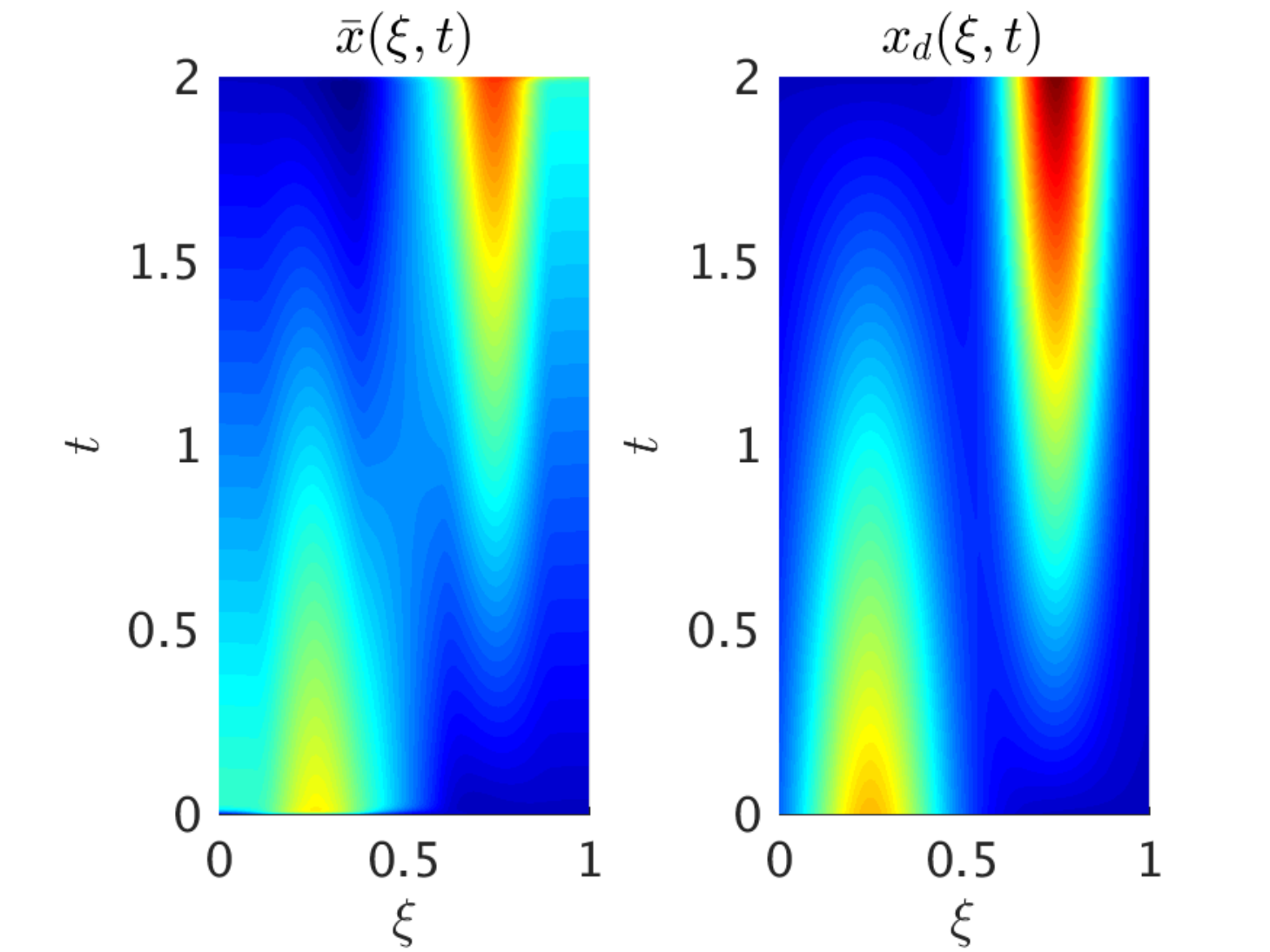}
			\includegraphics[width=0.49\linewidth]{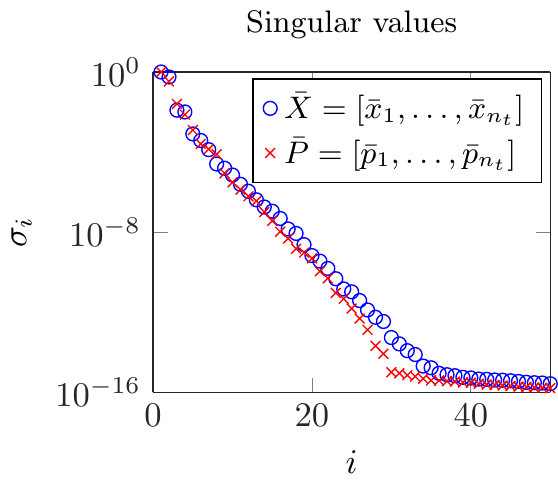}
			\caption{1D heat equation with $n=1000$ and $n_t=2000$. \textit{Left}. Controlled and desired states. \textit{Right}. Singular values of solutions. }
			\label{fig:motivation_1d_heat}
		\end{center}
	\end{figure}
	Figure~\ref{fig:motivation_1d_heat} (right) shows a rapid (exponential) decay of the singular values of the state and adjoint matrices, indicating that a low-rank approximation of the form \eqref{eq:low_rank_decomp} can be computed for any desired accuracy.
\end{example}

\setlength\fheight{0.24\linewidth}
\setlength\fwidth{0.36\linewidth}

We now may utilize the approximations \eqref{eq:low_rank_decomp} within our favourite and most suitable Krylov-subspace solver. In the case of linear PDE-constraints this method performs very well. As the Riccati equation is nonlinear one needs to investigate how the procedure changes for nonlinear equations. Our approach here is motivated by recent work in \cite{DolS15}, where one performs a nonlinear iteration as an outer solver, e.g., a Newton or Picard iteration, and then solves a linearized space-time discretized equation in low-rank form. Of course the structure of the linearized systems is more complicated than the one given above. Typically, the nonlinearity adds more terms to the sum of the discretized matrices. The number of terms is depending on the rank of the solution from the previous step of the nonlinear iteration. We derive the structure of the linearized equations for our problem in Section \ref{subsec::LRForm}. The goal of our solver is to decouple the spatial and temporal degrees of freedom since every increase in $n$ and $n_t$ will result in a dramatic increase of the computational cost. To break the \textit{curse-of-dimensionality} we will use an outer non-linear solver and a space-time inner solver that is based on the tensor train format introduced in Section \ref{sec::tt}.

\section{A low-rank all-at-once method for the differential Riccati equation}
\label{sec::lrdre}

In this section, we extend the discussion from Section \ref{sec::lrol} on open-loop control to the case of the closed-loop variant of problem \eqref{jyu2}. While the open-loop control problem requires the solution of a linear system in saddle point form we here face the more challenging differential Riccati equations.
We begin with a summary of the theoretical foundation of optimal closed-loop control via differential Riccati equations. For a numerical example from PDE boundary control, we visualize the approximability of the time-varying solution $P(\cdot)$ of \eqref{eq:dre} in terms of the singular values of corresponding matrix unfoldings. Based on this observation, we derive a novel tensor-based Newton-Kleinman method.

\subsection{Optimal feedback control: low-rank or not low-rank}

Let us again consider an optimal control problem of the form 
\begin{equation}
\label{jyu2_cl}
\begin{aligned}
\min_{\substack{x\in L^2(0,t_f;\mathbb R^n)\\u\in L^2(0,t_f;\mathbb R^m)}}
J(x,u):&=\frac{1}{2}\int\limits_0^{t_f} \| y(t) - Cx_d(t) \|^2 \, \mathrm{d}t+\frac{\beta}{2}\int\limits_0^{t_f}\| u(t)\|^2 \,  \mathrm{d}t + \frac{1}{2} \| y(t_f)- Cx_d(t_f) \|^2 \\
\text{s.t. } \  E\dot{x}(t)&= Ax(t) + Bu(t), \ \ x(0)=x_0.
\end{aligned}
\end{equation}
With the intention of computing controls that are robust (e.g., with respect to perturbations of the initial state), we are interested in a feedback control solving \eqref{jyu2_cl}.  It is well-known, see, e.g., \cite[Chapter 5]{Hin02} that the optimal control $u_{\mathrm{opt}}$ in fact can be expressed as 
\begin{align*}
u_{\mathrm{opt}}(t) = - B^T PE (x(t)-x_d(t)) - B^T r,
\end{align*}
where $P$ and $r$ satisfy
\begin{align}
E^T\dot{P}E + A^T PE + E^T PA  -E^T PBB^TPE + C^T C &= 0, \tag{DRE} \ \ 
E^TP(t_f)E = C^T C, \label{dre} \\
E^T \dot{r} + (A-BB^TPE)^T r + E^T P(Ax_d - E\dot{x}_d) &= 0, \notag \ \ r(t_f) = 0.
\end{align}
Note that both of the previous equations are \emph{initialized} at time $t_f$ and run backwards in time. Since the computation of the unknown $r\colon [0,t_f]\to \mathbb R^n$ can be realized by the low-rank technique from \cite{StoB15}, our primary focus is on the differential Riccati equation. For an efficient tensor based approach, the problem of interest has to exhibit a certain low-rank property. With that in mind, below we present a prototypical PDE-constrained optimal control problem and investigate the ranks of its low-rank tensor decomposition based in the tensor train format introduced later.

\begin{figure}[htb]
	\setlength\fwidth{0.20\linewidth} 
	\setlength\fheight{0.20\linewidth}
		\includegraphics[scale=0.9]{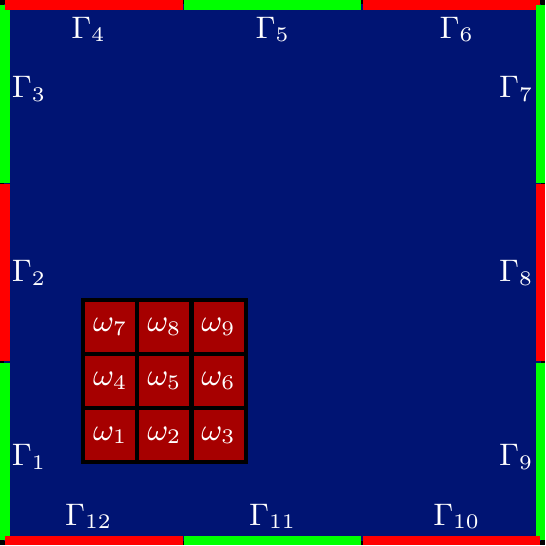}
	\setlength\fwidth{0.395\linewidth} 
	\setlength\fheight{0.395\linewidth}
	\includegraphics[scale=0.9]{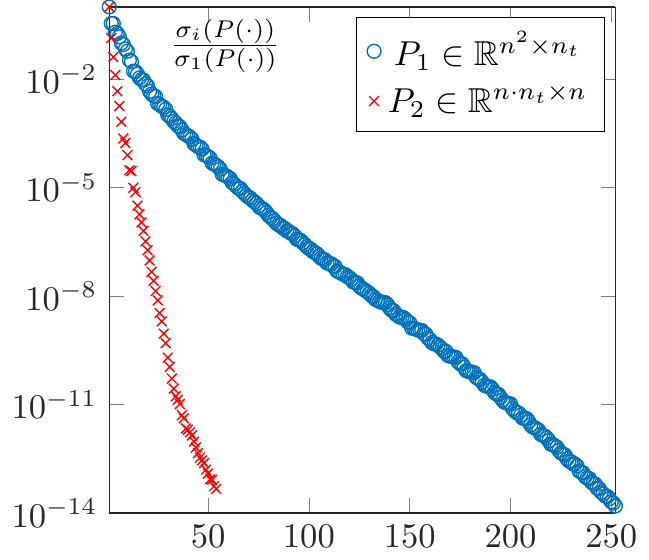}
	\caption{\emph{Left}. Control and observation domains. \emph{Right}. Singular values of matrix unfoldings.}
	\label{fig:cont_patches}
\end{figure}

\begin{example}
	Consider the two-dimensional heat equation subject to a Neumann boundary control of the form
	\begin{align*} 
	\frac{\partial}{\partial t} x(\xi,t) &= \Delta_{\xi} x(\xi,t)  && \text{in } \Omega \times (0,t_f), \\
	\frac{\partial}{\partial \nu} x(\xi,t)&=0=\sum_{i=1}^m \alpha_i(\xi)u(t) && \text{on } \Gamma \times (0,t_f), \\
	x(\xi,0)&=x_0(\xi) && \text{in } \Omega,
	\end{align*}
	where
	$\Omega=(0,1)^2,\Gamma=\partial \Omega$ and the control patches $\alpha_i,i=1,\dots,12$ are piecewise constants that are located on $\Gamma$, see Figure \ref{fig:cont_patches}. We further define a cost functional of the form 
	\begin{align*}
	\min_{\substack{u\in L^2(0,t_f;\mathbb R^{12})}}
	J(x,u):&=\frac{1}{2}\int\limits_0^{t_f} \| y(t) - e \|_{\mathbb R^9}^2 \, \mathrm{d}t+\frac{\beta}{2}\int\limits_0^{t_f}\| u(t)\|_{\mathbb R^{12}}^2 \,  \mathrm{d}t.
	\end{align*}
	with $e=(1,\dots,1)^T$ and $y_i(t)=\frac{1}{|\omega_i|} \int _{\omega_i} x(\cdot,t)\  \mathrm{d}t, i=1,\dots,9$. We have implemented a spatial discretization of the problem with $n=1089$ piecewise linear finite elements and a temporal implicit Euler scheme with $n_t=1000$ points. In Figure \ref{fig:cont_patches}, we show the singular values of the matrix unfoldings 
	\begin{align*}
	P_1=[\mathrm{vec}(P(t_1)),\dots,\mathrm{vec}(P(t_{n_t}))] \in \mathbb R^{n^2\times n_t}, \quad P_2=\begin{bmatrix} P(t_1) \\ \vdots \\ P(t_{n_t}) \end{bmatrix} \in \mathbb R^{n\cdot n_t \times n}.
	\end{align*} 
	In both cases, the numerical ranks $r(P_1)=252$ and $r(P_2)=54$ (for machine precision) are significantly smaller than the maximum rank of $r=1000$, allowing us to approximate the tensor $P\in \mathbb R^{n\times n\times n_t}$ by a low-rank representation $\tilde{P}$. Some results for such an approximation are given in Figure \ref{fig:storage_red} where a full representation of the solution $P$ is compared to a tensor train approximation $\tilde{P}$ (cf. Section \ref{sec::tt}).
\end{example}

\begin{figure}[htb]
	\begin{center}
		\begin{tabular}{|c|c|c|c|c|}
			\hline
			$\frac{\|P-\tilde{P}\|}{\| P\|}$ & $r(P_1)$ & $r(P_2)$ & $\#\texttt{kB}(\tilde{P})$ & $\frac{\#\texttt{kB}({P})}{\#\texttt{kB}(\tilde{P})}$  \\ \hline
			$\varepsilon$ & 252 & 54 & $1.21\cdot10^{5}$ & 78 \\
			10$^{-14}$ & 250 & 52 & $1.16\cdot10^{5}$ &  82  \\
			10$^{-12}$ & 216 & 37 & $7.18\cdot10^4$ & 132 \\
			10$^{-10}$ & 176 & 29 & $4.62\cdot10^{4}$ & 205  \\ 
			10$^{-8}$ & 133 & 22 & $2.70\cdot10^{4}$ & 354  \\
			10$^{-6}$ & 85 & 15 & $1.20\cdot10^{4}$ &  793 \\
			10$^{-4}$ & 45 & 8 & $3.59\cdot 10^3$ & 2642  \\
			10$^{-2}$ & 19 & 3 & $6.87\cdot 10^2$ & 13810  \\
			\hline
		\end{tabular}
		\caption{Storage reduction by tensor truncation from $P(\cdot)$ to $\tilde{P}(\cdot)$}
		\label{fig:storage_red}
	\end{center}
\end{figure}

\subsection{A space-time Newton-Kleinman formulation}
\label{subsec::LRForm}

Based on the \- well-known idea of a Newton-Kleinman approach, in the following we derive a nonlinear iteration for the fully discrete version of the differential Riccati equation.  At the $i$th Newton-Kleinman iteration ($i=0,1,\dots$), we need to solve
\begin{align*}
E^T \dot{P}_{i+1}E + \mathcal{A}(P_i)^T P_{i+1} E + E^TP_{i+1} \mathcal{A}(P_i) +C^TC + E^TP_i BB^T P_iE &= 0,
\\   E^T P_{i+1}(t_f)E &= C^T C,
\end{align*}
for $P_{i+1}$, where $\mathcal{A}(P_i)=A-BB^TP_iE$ denotes the closed-loop system operator associated with the current iteration. Note that each step still requires solving a differential matrix equation for the time-varying unknown $P_{i+1}\colon [0,t] \to \mathbb R^{n\times n}.$ However, in contrast to \eqref{dre} the resulting \emph{differential Lyapunov equation} is linear. Applying the vectorization operator to both sides of the above equation, we obtain the equivalent equation
\begin{align*}
\mathcal{E}   \dot{\mathcal{P}}_{i+1} + 
(\mathcal{L}+\mathcal{M}(\mathcal{P}_i)) \mathcal{P}_{i+1} +  \mathcal{C}+ 
\mathcal{G}(\mathcal{P}_i) = 0, \ \mathcal{E}\mathcal{P}_{i+1}(t_f) = 
\mathcal{C}
\end{align*}
on $\mathcal{P}_{i+1} = \mathrm{vec}(P_{i+1})$, where
\begin{align*}
\mathcal{E}&= E^T \otimes E^T, \quad  \mathcal{L} = E^T\otimes A^T + A^T \otimes 
E^T, \\[0.5ex]
\mathcal{M}(\mathcal{P}_i) &= E^T\otimes E^TP_i BB^T + E^TP_i BB^T 
\otimes 
E^T \\[0.5ex]
\mathcal{C} &= \mathrm{vec}(C^TC), \ \ \mathcal{G}(\mathcal{P}_i) = 
(E^TP_i \otimes E^T P_i) \mathcal{B}, \ \ 
\mathcal{B}=\mathrm{vec}(BB^T).
\end{align*}
With the intention of an all-at-once approach, we continue with a time stepping scheme for the interval $[0,t_f]$ with $n_t+1$ grid points and equidistant grid size $\tau=\frac{t_f}{n_t}$. Let us exemplarily consider an implicit Euler scheme leading to an inner iteration of the form
\begin{align*}
\frac{1}{\tau}\mathcal{E}(\mathcal{P}_{i+1}^{j+1}-\mathcal{P}_{i+1}^j)+ 
(\mathcal{L}+\mathcal{M}(\mathcal{P}_i^j)) \mathcal{P}_{i+1}^j +  \mathcal{C}+ 
\mathcal{G}(\mathcal{P}_i^j) &= 0, \ \ j =1,\dots,n_t, \\
\mathcal{E}\mathcal{P}_{i+1}^{n_t+1} &= \mathcal{C}. 
\end{align*}
Let us emphasize that the subscript $i$ indicates the iteration of the outer non-linear solver and the superscript $j$ enumerates the time-step of the current vector.
Note that the above scheme indeed is implicit since the equation is running backwards in time. Similar as before, we represent all time steps in a single vector that we denote by 
\begin{align*}
\mathbf{p}_i : = \begin{bmatrix} \mathcal{P}_i^1 \\ \vdots \\ \mathcal{P}_i^{n_t} \end{bmatrix}\in\R^{n^2n_t}. 
\end{align*}
As a consequence, we can express the inner iteration as a tensorized linear system
\begin{align*}
(\mathbf{L}+ \mathbf{D}+\mathbf{M}(\mathbf{p}_i) )\mathbf{p}_{i+1} = 
\mathbf{c} + \mathbf{g}(\mathbf{p}_i),
\end{align*}
where
\begin{align*}
\mathbf{L} &= I_{n_t}\otimes \mathcal{L}, \ \  \mathbf{D}=D\otimes 
\mathcal{E}, \ \ \mathbf{M}(\mathbf{p}_i)  = \mathrm{blkdiag}( 
\mathcal{M}(\mathcal{P}_i^1),\ldots,\mathcal{M}(\mathcal{P}_i^{n_t}) )  \\[2ex]
D&=\frac{1}{\tau}\begin{bmatrix} -1 & 1 &  & \\ &\ddots & \ddots & \\ & & 
-1 & 1 \\ & & & -1 \end{bmatrix}, \ \ \mathbf{c}=-\begin{bmatrix} \mathcal{C} 
\\ \vdots \\ \mathcal{C} \\ \mathcal{C}+\frac{1}{\tau} \mathcal{C} 
\end{bmatrix}, \ \ \mathbf{g}(\mathbf{p}_i) = - \begin{bmatrix} 
\mathcal{G}(\mathcal{P}_i^1) \\ \vdots \\ \vdots \\  
\mathcal{G}(\mathcal{P}_i^{n_t}) \end{bmatrix}.
\end{align*}
For the main idea, we may think of $\mathbf{p}_i$ as being a tensor which we approximate. For this we realize that we can write 
$$
\mathcal{P}_i^j=\left(e_j^T\otimes I_n\otimes I_n\right)\mathbf{p}_i.
$$
Let us first assume that
$$
\mathbf{p}_i  \approx 
u \otimes v \otimes  w,
$$ 
where $u \in \mathbb{R}^{n_t}$ and $v,w \in \mathbb{R}^{n}$ are column vectors.
We then get
$$
\mathcal{P}_i^j\approx \left(e_j^T\otimes I_n\otimes I_n\right)\left(u \otimes v \otimes  w\right)=
u_j\left( v \otimes  w\right).
$$
For obtaining a matrix valued expression, let us interpret $\mathcal{P}_i^j=\mathrm{vec}(P_i^j)$ as the vectorization of a matrix. It then follows that
$$
\mathrm{vec}(P_i^j)\approx u_j (v\otimes w)\mathrm{vec}(1)\ \Leftrightarrow  \ \ P_i^j \approx u_j wv^T.
$$
We use the latter expression for a derivation of the term $\mathcal{M}(\mathcal{P}_i^j)$ which is of the form
\begin{align*}
\mathcal{M}(\mathcal{P}_i^j) &=  E^T\otimes E^TP_i^j BB^T + E^TP_i^j BB^T \otimes E^T \\
&\approx  E^T\otimes E^T\left(u_jwv^T\right) BB^T + E^T\left(u_j wv^T\right) BB^T \otimes E^T \\
&\approx u_j \left[E^T\otimes E^T\left( wv^T\right) BB^T + E^T\left(wv^T\right) BB^T \otimes E^T \right].
\end{align*}
With this in mind, let us now move beyond the case of a simple rank-one approximation and utilize the tensor train format illustrated in Section \ref{sec::tt}. This format allows the approximation of the state tensor via
%
$$
\mathbf{p}_i= \!\!\sum_{s_1,s_2=1}^{r_1,r_2} 
\mathbf{u}^{(1)}_{s_1} \otimes \mathbf{u}^{(2)}_{s_1,s_2} \otimes \mathbf{u}^{(3)}_{s_2}
$$
with $r_1$ and $r_2$ the ranks of the tensor train approximation.
As a result, we can approximate the operator $\mathbf{M}$ via 
\begin{align*}
\mathbf{M}(\mathbf{p}_i) &\approx \sum_{s_1,s_2=1}^{r_1,r_2} 
\underset{k=1,\dots,n_t}{\mathrm{blkdiag}}(  
E^T\otimes E^Tu_{s_1,k}^{(1)} \mathbf{u}^{(3)}_{s_2}\left(\mathbf{u}^{(2)}_{s_1,s_2}\right)^T BB^T + E^Tu_{s_1,k}^{(1)} \mathbf{u}^{(3)}_{s_2}\left(\mathbf{u}^{(2)}_{s_1,s_2}\right)^T BB^T \otimes E^T  
) \\
&= \sum_{s_1,s_2=1}^{r_1,r_2}  \underbrace{\mathrm{diag}(u_{s_1}^{(1)} )}_{=:U_j} \otimes ( E^T \otimes E^T 
\mathbf{u}^{(3)}_{s_2}\left(\mathbf{u}^{(2)}_{s_1,s_2}\right)^T BB^T +  E^T \mathbf{u}^{(3)}_{s_2}\left(\mathbf{u}^{(2)}_{s_1,s_2}\right)^T BB^T \otimes E^T  ) \\[1ex]
\intertext{Similarly, using the \emph{Hadamard product} of two vectors, we obtain}
\mathbf{g}(\mathbf{p}_i) &\approx -  \left[
\left(
\sum_{j=1}^{r_1} E^T u_{j,k}^{(1)} \mathbf{u}^{(3)}_{1:r_2} \left(\mathbf{u}^{(2)}_{j,1:r_2}\right)^T
\otimes \sum_{\ell=1}^{r_1} E^T  u_{l,k}^{(1)} \mathbf{u}^{(3)}_{1:r_2} \left(\mathbf{u}^{(2)}_{l,1:r_2}\right)^T
\right)\mathcal{B}\right]_{k=1}^{n_t} \\
&= -\left[\left(E^T \mathbf{u}^{(3)}_{1:r_2} \otimes E^T \mathbf{u}^{(3)}_{1:r_2} \right) \mathrm{vec}\left(\sum_{j,\ell=1}^{r_1} \underbrace{(u_{j,k}^{(1)} \circ u_{\ell,k}^{(1)} )}_{=:(u_{j\circ \ell})_k} \left(\mathbf{u}^{(2)}_{j,1:r_2}\right)^T BB^T \mathbf{u}^{(2)}_{\ell,1:r_2} \right)\right]_{k=1}^{n_t} \\
&= -\sum_{j,\ell=1}^{r_1} u_{j\circ \ell} \otimes \left(E^T \mathbf{u}^{(3)}_{1:r_2} \otimes E^T \mathbf{u}^{(3)}_{1:r_2}\right) \mathrm{vec}\left(\left(\mathbf{u}^{(2)}_{j,1:r_2}\right)^T BB^T \mathbf{u}^{(2)}_{\ell,1:r_2} \right).
\end{align*}
Note that we have included only one sum in the derivation of $\mathbf{g}(\mathbf{p}_i)$ and have collected the sum over $s_2$ into the matrix products with $ \mathbf{u}^{(3)}_{1:r_2} \left(\mathbf{u}^{(2)}_{j,1:r_2}\right)^T,$ where $1:r_2$ indicates $r_2$ vectors collected into a $n \times r_2$ matrix, which is standard Matlab notation.
Altogether, we obtain a linear system of tensor product structure of the form 
$\mathbf{A}\mathbf{p}_{i+1} = \mathbf{f},$ where 
\begin{align}
\label{system1}
\mathbf{A} &= I_{n_t} \otimes E^T \otimes A^T + I_{n_t}\otimes A^T \otimes E^T
+ D \otimes E^T \otimes E^T           +\mathbf{M}(\mathbf{p}_i) , \\[2ex]
\mathbf{f} &= -\sum_{k=1}^{r_f} \mathbf{e}_{\tau} \otimes  c_k^T \otimes c_k^T +\mathbf{g}(\mathbf{p}_i)
\end{align}
where $\mathbf{e}_\tau = [1,\dots,1,1+\frac{1}{\tau}]^T \in \mathbb R^{n_t}.$ 		
The resulting system is now solved with the Alternating Minimal Energy (AMEn) solver, presented next. Note that the number of terms in $\mathbf{M}(\mathbf{p}_i)$ depends on the previous solution and we will illustrate in the numerical experiments how this behaves with respect to changes in the system parameters and system size. We will use a tolerance $\varepsilon_{\mathrm{newton}}$ to assess the   convergence of the Newton iterates via 
$
\frac{\norm{\mathbf{p}_{i+1}-\mathbf{p}_{i}}}{\norm{\mathbf{p}_{i+1}}}< \varepsilon_{\mathrm{newton}}.
$
\subsection{Nested approach}
As the nonlinear iteration can become quite expensive with a possibly high rank of the state variables we aim to drastically reduce this cost by using a nested approach \cite{HerS10,PeaSW14}. There the authors solve a nonlinear problem on a coarse mesh and then transfer the solution to the next finer mesh as the initial guess for the nonlinear iteration. This process is continued until the final mesh-size is solved. Typically, this reduces the number of nonlinear iterations drastically.
Our aim is to apply this approach here as well. We will solve the low-rank DRE on a coarse mesh to obtain a solution that we then transfer to the next finer mesh as the intial guess for the Newton--Kleinman iteration.

The crucial ingredient in such a multilevel strategy is the prolongation or interpolation of the coarse solution onto the finer grid. In \cite{grasedyck2007multigrid,vandereycken2009local} the authors discuss such a multilevel strategy for large scale matrix equations. For the combination of multigrid with tensor methods we refer to \cite{khoromskij2009multigrid}. 

We assume that the univariate prolongation is performed via the multiplication with a matrix $\mathbf{\Pi}_h\in\R^{n_{h},n_{2h}},$ where
$n_{h},n_{2h}$ are the dimensions of the fine and coarse mesh, respectively. Then the prolongation that is required in the matrix equation setup is defined via
$$
\mathbf{\Pi}=\mathbf{I}_{n_t} \otimes \mathbf{\Pi}_h\otimes \mathbf{\Pi}_h
$$
and if we have computed a solution on grid level $2h$ as
\begin{equation}
\mathbf{p}_{2h}=\!\!\sum_{s_1,s_2=1}^{r_1,r_2} 
\mathbf{p}^{(1)}_{s_1} \otimes \mathbf{p}^{(2)}_{s_1,s_2} \otimes \mathbf{p}^{(3)}_{s_2}
\end{equation}
we then obtain 
\begin{equation}
\mathbf{p}_{h}=\mathbf{\Pi}\mathbf{p}_{2h}=\!\!\sum_{s_1,s_2=1}^{r_1,r_2} 
\mathbf{p}^{(1)}_{s_1} \otimes \mathbf{\Pi}_h\mathbf{p}^{(2)}_{s_1,s_2} \otimes \mathbf{\Pi}_h\mathbf{p}^{(3)}_{s_2}.
\end{equation}
This will now be a low-rank solution on the fine grid and it will be the starting guess for the Newton--Kleinman iteration on the finer mesh. It is obvious that this process can be continued further.

\section{The tensor train decomposition and algorithms}
\label{sec::tt}
The efficient solution of tensor-valued equations has recently seen much progress \cite{Tob12,BalG13,Ose11b,GraKT13,KolB09,hackbusch-2012}. While there is a variety of available tensor formats we decide to use the so-called tensor train (TT) representation \cite{Ose11,Ose11b,OseT09}. As our problem essentially consists of tensors of order three, the Tucker format \cite{GraKT13,KolB09} would also be appropriate. The availability of tailored solvers along with the possibility of designing and using preconditioners makes the TT-format and its toolbox\footnote{\url{https://github.com/oseledets/TT-Toolbox}} \cite{tt-toolbox} an ideal candidate for our purposes. 

In the following we will briefly introduce the TT format, while we point to the literature for details.  To this end, suppose $\mathbf{p}\in\R^{n^2 n_t}$ is the approximate solution vector of a space-time matrix equation as discussed in Section \ref{sec::lrdre}.
However, its elements can be also naturally enumerated by three indices $i_1,i_2,i_3,$ corresponding to the discretization in time and the two spatial dimensions, respectively.
Introducing a \emph{multi-index} 
$$
\overline{i_1i_2i_3} = (i_1-1) n^2 +(i_2-1)n + i_3,
$$
we can denote $\mathbf{p} = \left[\mathbf{p}(\overline{i_1i_2i_3})\right]_{i_1,i_2,i_3=1}^{n_t,n,n}$, and consider $\mathbf{p}$ as a three-dimensional \emph{tensor} with elements $\mathbf{p}(i_1,i_2,i_3)$.
The TT decomposition aims to represent $\mathbf{p}$ as follows
\begin{equation}
\label{eq:tt}
\mathbf{p}(i_1,i_2,i_3) = \!\!\sum_{s_1,s_2=1}^{r_1,r_2} \mathbf{u}^{(1)}_{s_1}(i_1) \mathbf{u}^{(2)}_{s_1,s_2}(i_2) \mathbf{u}^{(3)}_{s_2}(i_3) \,\, \Leftrightarrow \,\,
\mathbf{p} = \!\!\sum_{s_1,s_2=1}^{r_1,r_2}
\mathbf{u}^{(1)}_{s_1} \otimes \mathbf{u}^{(2)}_{s_1,s_2} \otimes \mathbf{u}^{(3)}_{s_2}.
\end{equation}

The ranges $r_1,r_2$ are called TT ranks and the $\mathbf{u}^{(m)}$, $m=1,2,3$
are the so-called TT blocks, {\color{black}with}
$\mathbf{u}^{(1)}\in\mathbb{R}^{n_t \times r_1}$, $\mathbf{u}^{(2)}\in\mathbb{R}^{r_1 \times n \times r_2}$ and $\mathbf{u}^{(3)}\in\mathbb{R}^{r_2 \times n}$.
When fixing the indices we get for $\mathbf{u}^{(2)}(i_2)\in\mathbb{R}^{r_1 \times
	r_2}$ a matrix slice, for $\mathbf{u}^{(2)}_{s_1,s_2}\in\mathbb{R}^{n}$ a 
vector, and for $\mathbf{u}^{(2)}_{s_1,s_2}(i_2)$ a scalar. The values of
$r_1,r_2$ depend on the accuracy enforced in equation 
\eqref{eq:tt}.
Given a full tensor or a TT decomposition with excessively large ranks, the quasi-optimal TT approximation can be computed using $d-1$ Singular Value Decompositions (SVD), truncated up to the desired accuracy threshold~\cite{Ose11}.
For example, we can reduce $r_1$ by computing SVD $\mathbf{u}^{(1)} = U \Sigma V^T$, and select only the first $r_1'<r_1$ singular vectors, approximating $\mathbf{\tilde u}^{(1)} = U_{1:r_1'}$ and $\mathbf{\tilde u}^{(2)}(i_2) = \Sigma_{1:r_1'} V_{1:r_1'}^T \mathbf{u}^{(2)}(i_2)$.
We introduce a \emph{truncation tolerance} $\varepsilon_{\mathrm{trunc}}>0$ and select $r_1'$ such that $\|\Sigma_{1:r_1'} - \Sigma\|_F \le \varepsilon_{\mathrm{trunc}} \|\Sigma\|_F$.
The other TT ranks can be truncated similarly.

The same technique can be used to represent a linear operator $\mathbf{A}$ in TT format,
\begin{equation}\label{eq:ttm}
\mathbf{A} \approx \sum_{\ell_1,\ell_2=1}^{R_1,R_2} A_{\ell_1}^{(1)} \otimes A_{\ell_1,\ell_2}^{(2)} \otimes A_{\ell_2}^{(3)}.
\end{equation}
Note that the Kronecker product matrix assembly \eqref{system1} is a particular case of \eqref{eq:ttm}.

For the numerical solution of $\mathbf{A}\mathbf{p}=\mathbf{f}$
we use the Alternating Minimal Energy (AMEn) algorithm~\cite{DolS14}, which is an enhanced version of the Alternating Linear Scheme (ALS)~\cite{HolRS12}.
The initial idea is to rewrite the TT decomposition~\eqref{eq:tt} as a linear map from the elements of a single TT block, and consider $\mathbf{A}\mathbf{p}=\mathbf{f}$ as an overdetermined system on these elements.
Specifically for the three-dimensional TT format, we introduce \emph{frame} matrices
\begin{align}
U_{\neq 1} & = \sum_{s_2=1}^{r_2}
I_{n_t} \otimes \left(\mathbf{u}^{(2)}_{1:r_1,s_2}\right)^T \otimes \mathbf{u}^{(3)}_{s_2} & \in \mathbb{R}^{n^2n_t \times n_t r_1}, \\
U_{\neq 2} & = \mathbf{u}^{(1)} \otimes I_n \otimes \left(\mathbf{u}^{(3)}\right)^T & \in \mathbb{R}^{n^2 n_t \times r_1 n r_2}, \\
U_{\neq 3} & = \sum_{s_1=1}^{r_1}
\mathbf{u}^{(1)}_{s_1} \otimes \mathbf{u}^{(2)}_{s_1,1:r_2} \otimes I_n & \in\mathbb{R}^{n^2n_t \times r_2n},
\end{align}
as well as the vectorisations of the trailing TT blocks $\mathrm{vec}_T(\mathbf{u}^{(1)}):=\mathrm{vec}((\mathbf{u}^{(1)})^T)$, $\mathrm{vec}_T(\mathbf{u}^{(3)}):=\mathrm{vec}((\mathbf{u}^{(3)})^T)$ with $\mathrm{vec}(\cdot)$ defined in Sec.~1, and of the middle TT block,
$$
\mathrm{vec}_T(\mathbf{u}^{(2)}) := \mathrm{vec}\begin{bmatrix}(\mathbf{u}^{(2)}(1))^T\\ \vdots \\ (\mathbf{u}^{(2)}(n))^T\end{bmatrix}\in \mathbb{R}^{r_1nr_2}.
$$
One can notice that, given~\eqref{eq:tt}, the identity $\mathbf{p} = U_{\neq k} \mathrm{vec}_T(\mathbf{u}^{(k)})$ holds for any $k=1,2,3$.
Now we can plug this \emph{reduced basis} ansatz into the original system $\mathbf{A}\mathbf{p}=\mathbf{f}$, and resolve it by \emph{projecting} onto the same frame matrix,
\begin{equation}\label{eq:localsys}
(U_{\neq k}^T \mathbf{A} U_{\neq k}) \mathrm{vec}_T(\mathbf{u}^{(k)}) = U_{\neq k}^T \mathbf{f}.
\end{equation}
Now, the ALS algorithm iterates over $k=1,2,3$, solving~\eqref{eq:localsys} in each step, and updating the TT block $\mathbf{u}^{(k)}$.
This simple algorithm requires an initial guess in the TT format~\eqref{eq:tt} with fixed TT ranks, which might be difficult to guess a priori.
To circumvent this issue, AMEn algorithm computes additionally a TT approximation of the current residual $\mathbf{f}-\mathbf{A}\mathbf{p}$ in a form similar to~\eqref{eq:tt} with (smallish) TT ranks $\rho_1,\rho_2$, and expands the TT blocks of the solution with the TT blocks of this residual approximation~\cite{DolS14}.
This allows us to increase the TT ranks of the solution from $r_1,r_2$ to $r_1+\rho_1, r_2+\rho_2$ in each iteration and improve convergence.
The iteration continues until the accuracy and the TT ranks of the solution reach their desired values.
At the same time, we can descrease the ranks using SVD if they become too large.
Specifically, let $\mathbf{u}^{(1)},\mathbf{u}^{(2)},\mathbf{u}^{(3)}$ be the TT blocks at the given AMEn iteration, and let $\mathbf{\hat u}^{(1)},\mathbf{\hat u}^{(2)},\mathbf{\hat u}^{(3)}$ be the TT blocks after one full sweep over $k=1,2,3$.
We choose a \emph{stopping tolerance} $\varepsilon_{\mathrm{amen}}>0$ and stop the algorithm when
$$
\max\left[\frac{\|\mathbf{u}^{(1)} - \mathbf{\hat u}^{(1)}\|_F}{\|\mathbf{\hat u}^{(1)}\|_F},~\frac{\|\mathbf{u}^{(2)} - \mathbf{\hat u}^{(2)}\|_F}{\|\mathbf{\hat u}^{(2)}\|_F},~\frac{\|\mathbf{u}^{(3)} - \mathbf{\hat u}^{(3)}\|_F}{\|\mathbf{\hat u}^{(3)}\|_F}\right] \le \varepsilon_{\mathrm{amen}}.
$$
The same threshold $\varepsilon_{\mathrm{amen}}$ is used for SVD within the AMEn iteration when the TT ranks need to be reduced.

For the practical efficiency we can notice that the reduced matrix $B_k = U_{\neq k}^T \mathbf{A} U_{\neq k}$ in~\eqref{eq:localsys} can be written in a Kronecker product form that inherits the original matrix TT decomposition~\eqref{eq:ttm}.
In particular, we can write
$$
B_k=\sum_{\ell_1,\ell_2=1}^{R_1,R_2} B_{\ell_1}^{(1)}\otimes B_{\ell_1,\ell_2}^{(2)}\otimes B_{\ell_2}^{(3)},
$$
where $B^{(k)} = A^{(k)}$ (that corresponds to the identity factor in the frame matrix $U_{\neq k}$), and hence is large but sparse, while the other matrices are dense but small, of sizes $r_1\times r_1$, $r_2 \times r_2$ or even $1\times 1$ if $k=1$ or $k=3$ is considered.
This allows us to design an efficient block Jacobi preconditioner for a GMRES solver of~\eqref{eq:localsys}.
Without loss of generality, we can assume that $k=2$, in which case $B^{(2)}_{\ell_1,\ell_2} \in \mathbb{R}^{n \times n}$ is sparse, and $B^{(1)}_{\ell_1} \in \mathbb{R}^{r_1 \times r_1}$ and $B^{(3)}_{\ell_2} \in \mathbb{R}^{r_2 \times r_2}$ are dense.
Then we construct the block diagonal preconditioner as follows,
\begin{equation}\label{eq:localprec}
P_2 = \sum_{\ell_1,\ell_2=1}^{R_1,R_2} \mathrm{diag}(B_{\ell_1}^{(1)})\otimes B_{\ell_1,\ell_2}^{(2)}\otimes \mathrm{diag}(B_{\ell_2}^{(3)}).
\end{equation}
For $k=1$ or $k=3$ the procedure is similar.

\section{Numerical experiments}

In this section, we apply the above method to two control problems for the linear heat equation. The spatially discrete systems are obtained from a finite element discretization with piecewise linear elements and $n=289,1089,4225,16641$ degrees of freedom. The fully discrete system is obtained by an implicit Euler scheme with $n_t=1000$ equidistant time steps on the interval $[0,t_f]$, with $t_f=10$. All simulations are generated on an AMD Ryzen 7 1800X @ 3.68 GHz x 16, 64 GB RAM with  \matlab \;Version 9.2.0.538062 (R2017a). Our implementation is based on the tensor-train toolbox\footnote{\url{https://github.com/oseledets/TT-Toolbox}} and in particular we rely on the \texttt{amen\_block\_solve} function \cite{BenDOS15}, which allows the use of sparse matrices and the incorporation of the proposed preconditioner~\eqref{eq:localprec}. We set the truncation tolerance to $\varepsilon_{\mathrm{trunc}} = 10^{-12}$, and the stopping tolerance for AMEn to $\varepsilon_{\mathrm{amen}} = 10^{-9}$. The outer nonlinear Newton method is then solved up to a tolerance of $\varepsilon_{\mathrm{newton}} = 10^{-5}.$


\subsection{Incomplete distributed control}

\begin{figure}[htb]
	\begin{center}
		\setlength\fwidth{0.25\linewidth}
		\includegraphics[scale=0.8]{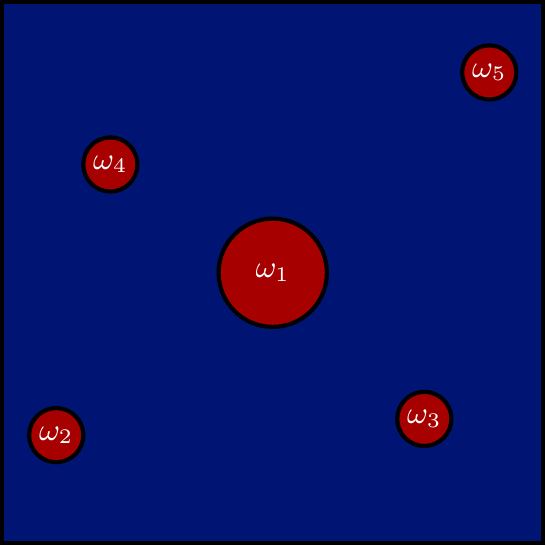}
		\setlength\fheight{0.30\linewidth}
		\includegraphics[scale=0.8]{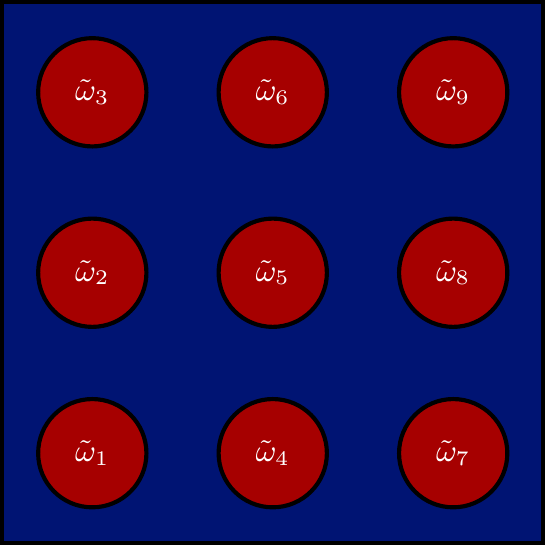}
		\caption{\emph{Left}. Control domains. \emph{Right}. Observation domains.}
		\label{fig:heat_cont_obsv_patches}
	\end{center}
\end{figure}

We consider the optimal control of the heat equation equipped with homogeneous Neumann boundary condition as given in
\begin{align*} 
\frac{\partial}{\partial t} x(\xi,t) &= \Delta_{\xi} x(\xi,t)+ \sum_{i=1}^5 \chi_{\omega_i}(\xi) u_i(t)  && \text{in } \Omega \times (0,t_f)  , \\
\frac{\partial}{\partial \nu} x(\xi,t)&=0 && \text{on } \Gamma \times (0,t_f), \\
x(\xi,0)&=x_0(\xi) && \text{in } \Omega.
\end{align*}
The spatial domain is given by $\Omega=(0,1)^2$ and we denote the boundary by $\Gamma=\partial \Omega$. As illustrated in Figure \ref{fig:heat_cont_obsv_patches} the control consists of five spatially constant control patches $\omega_1,\dots,\omega_5$ , given by 
\begin{align*}  
\omega_1 &= \{ \xi\in \Omega \ |\  \|\xi - \xi_{M_1} \|^2 \le 0.01 \}, \\
\omega_i &= \{ \xi \in \Omega \ |\  \|\xi - \xi_{M_i}  \|^2 \le 0.0025 \},\ i=2,3,4,5,
\end{align*}
with centres defined as 

$\xi_{M_1} =\begin{pmatrix} \frac{1}{2} \\[1ex] \frac{1}{2} \end{pmatrix}, 
\xi_{M_2} = \begin{pmatrix} 0.1 \\ 0.2 \end{pmatrix},  
\xi_{M_3} = \begin{pmatrix} 0.78 \\ 0.23 \end{pmatrix},  
\xi_{M_4} = \begin{pmatrix} 0.2 \\ 0.7  \end{pmatrix}, 
\xi_{M_5} = \begin{pmatrix} 0.9 \\ 0.87 \end{pmatrix}.$

\begin{figure}[h!]
	\begin{center}
		\includegraphics[scale=0.4]{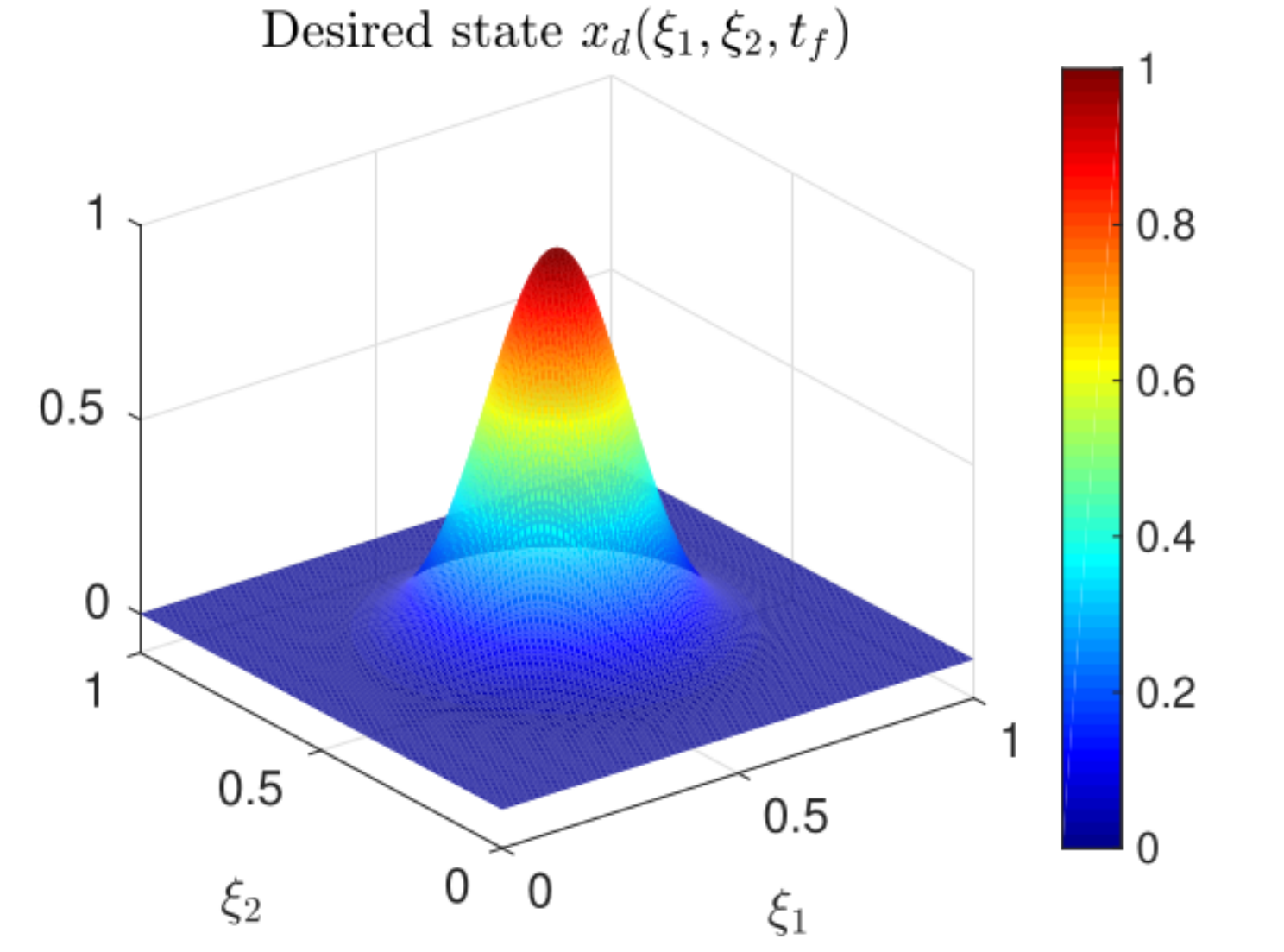}
		\includegraphics[scale=0.4]{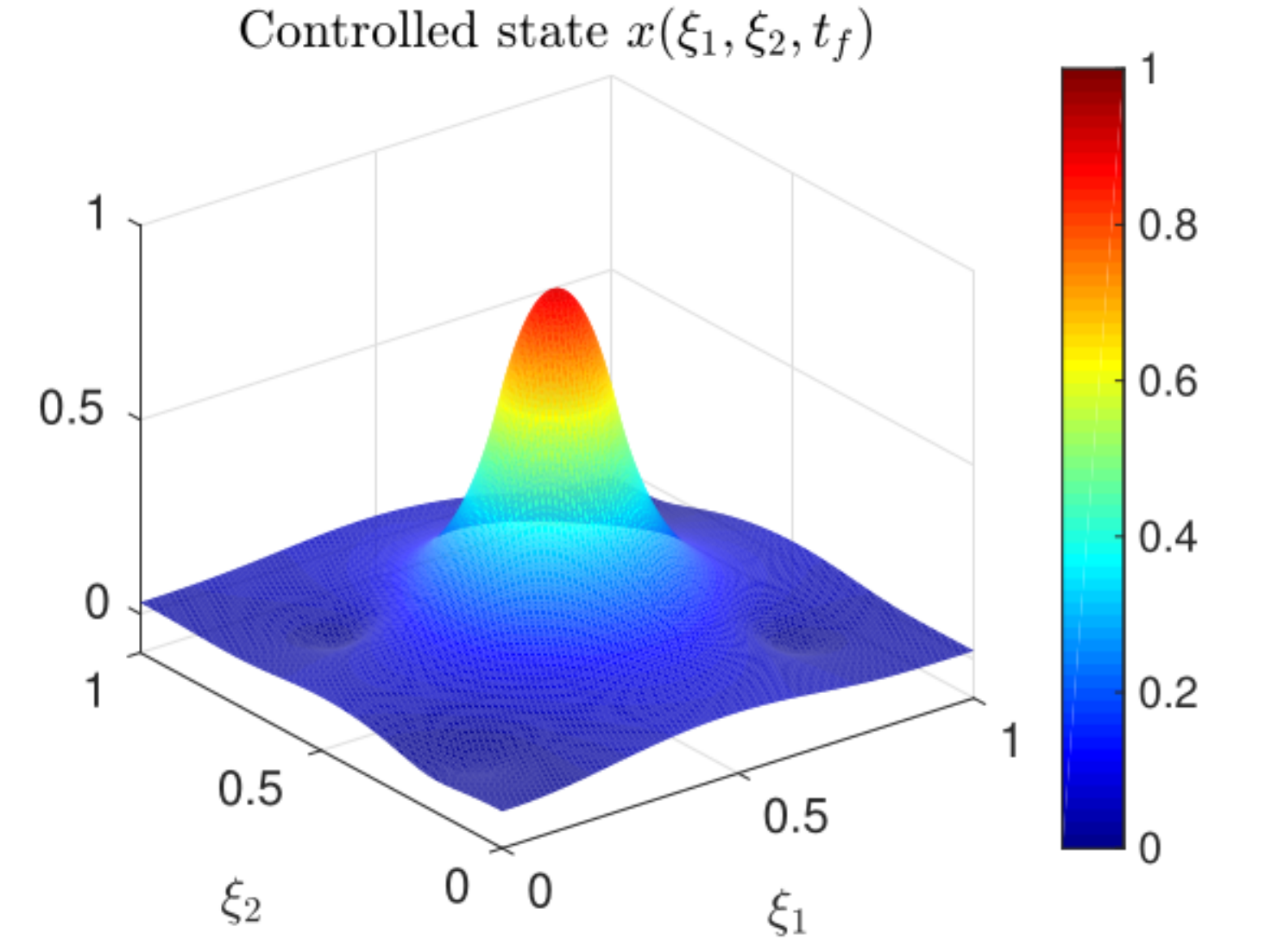}
		\caption{\emph{Left.} Desired state $x_d(\xi_1,\xi_2)$. \emph{Right.} Controlled state $x(\xi_1,\xi_2,t_f)$.}
		\label{fig:heat_int_des_vs_con_interior}
	\end{center} 
\end{figure}

The cost functional of interest is then given by 
\begin{align*}
\min_{\substack{u\in L^2(0,t_f;\mathbb R^{5})}}
J(x,u):&=\frac{1}{2}\int\limits_0^{t_f} \| y(t) - Cx_d \|_{\mathbb R^9}^2 \, \mathrm{d}t+\frac{\beta}{2}\int\limits_0^{t_f}\| u(t)\|_{\mathbb R^{5}}^2 \,  \mathrm{d}t,
\end{align*}
where $t_f= 10, \beta=10^{-4}.$ We depict the desired state $x_d$ in Figure \ref{fig:heat_int_des_vs_con_interior}, which we assume to be constant in time. For our problem the output $y(\cdot)$ is defined as
\begin{align*} 
y_i(t) &= \frac{1}{|\widetilde{\omega}_i|} \int_{\widetilde{\omega}_i} x(\xi,t) \, \mathrm{d}\xi , \text{where} \ \ 
\widetilde{\omega}_i = \{ \xi \in \Omega \ |\  \|\xi - \xi_{M_i} \|^2 \le 0.01 \}
\end{align*}
with the observational patches given via its centres
\begin{align*}
\xi_{M_1} &= \begin{pmatrix} \frac{1}{6} \\[1ex] \frac{1}{6} \end{pmatrix}, \
\xi_{M_2} = \begin{pmatrix} \frac{3}{6} \\[1ex] \frac{1}{6} \end{pmatrix}, \ 
\xi_{M_3} = \begin{pmatrix} \frac{5}{6} \\[1ex] \frac{1}{6} \end{pmatrix}, \
\xi_{M_4}= \begin{pmatrix} \frac{1}{6} \\[1ex] \frac{3}{6} \end{pmatrix}, \
\xi_{M_5} = \begin{pmatrix} \frac{3}{6} \\[1ex] \frac{3}{6} \end{pmatrix}, \\
\xi_{M_6} &= \begin{pmatrix} \frac{5}{6} \\[1ex] \frac{3}{6} \end{pmatrix}, \
\xi_{M_7} = \begin{pmatrix} \frac{1}{6} \\[1ex] \frac{5}{6} \end{pmatrix}, \  
\xi_{M_8} = \begin{pmatrix} \frac{3}{6} \\[1ex] \frac{5}{6} \end{pmatrix}, \
\xi_{M_9} = \begin{pmatrix} \frac{5}{6} \\[1ex] \frac{5}{6} \end{pmatrix}.
\end{align*}
These are illustrated in Figure \ref{fig:heat_cont_obsv_patches}. 

\begin{figure}[h!]
	\begin{center}
		\includegraphics[scale=1]{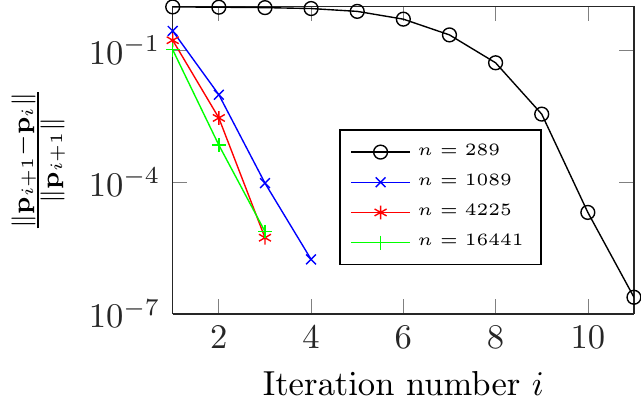}
		\includegraphics[scale=1]{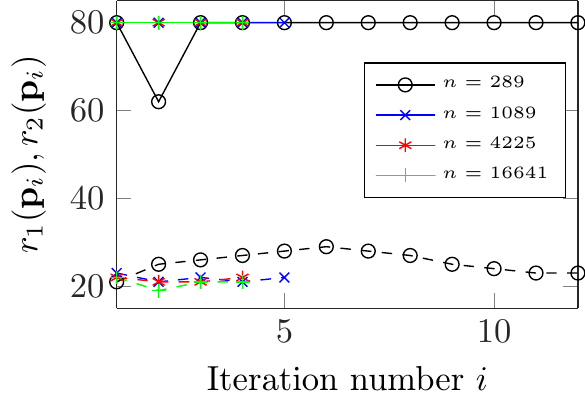}
		\caption{\emph{Left}. Change in Newton iterates. \emph{Right}. TT ranks of $\mathbf{p}_i$ during Newton iteration, $r_1$ (solid lines) and $r_2$ (dashed lines).}
		\label{fig:nk_iterates_ranks_interior}
	\end{center}
\end{figure}

In Figure \ref{fig:nk_iterates_ranks_interior} we show results for the all-at-once computation of the solution $\mathbf{p}\in \mathbb R^{n\times n \times 1000}$ with varying grid sizes $n=289,1089,4225,16641$. For the coarsest discretization, the Newton-Kleinman iteration is initialized with the zero solution, i.e., $\mathbf{p}_0=0$. 

As can be seen from Figure \ref{fig:nk_iterates_ranks_interior} the quadratic convergence of the Newton iterates is obtained. For finer discretization levels, the Newton-Kleinman iteration is combined with the nested approach discussed in the previous section. In contrast to the coarsest mesh case $n=289$, the number of outer Newton steps required to obtain the desired stopping criterion of $\varepsilon_{\mathrm{newton}} = 10^{-5}$ is significantly reduced.
Additionally, we also show the individual TT-ranks of the iterates throughout the iteration. Let us emphasize that the largest numerical rank $r=80$ is drastically smaller than the maximal possible ranks $r=1000$ and $r=16641$, respectively. In fact, storing the full solution tensor would require more than 2TB of data while the low-rank approximation only requires $\approx$230MB.
In Figure \ref{fig:heat_int_des_vs_con_interior}, we also show the results for the approximation of the final state associated with the solution of the optimal control problem. 

\newpage

\subsection{Boundary control}
\begin{figure}[htb]
	\begin{center}
\includegraphics[scale=1]{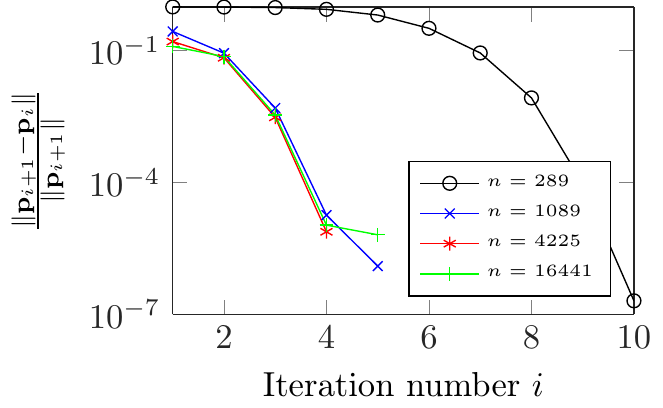}
		\includegraphics[scale=1]{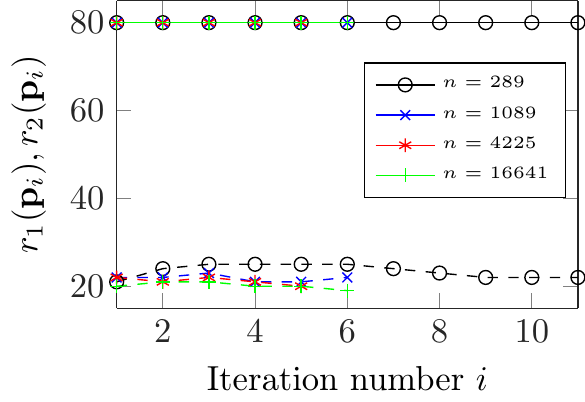}
		\caption{\emph{Left}. Change in Newton iterates. \emph{Right}. TT ranks of $\mathbf{p}_i$ during Newton iteration, $r_1$ (solid lines) and $r_2$ (dashed lines).}
		\label{fig:nk_iterates_ranks_boundary}
	\end{center}
\end{figure}

\begin{figure}[htb]
	\begin{center}
		\includegraphics[scale=0.4]{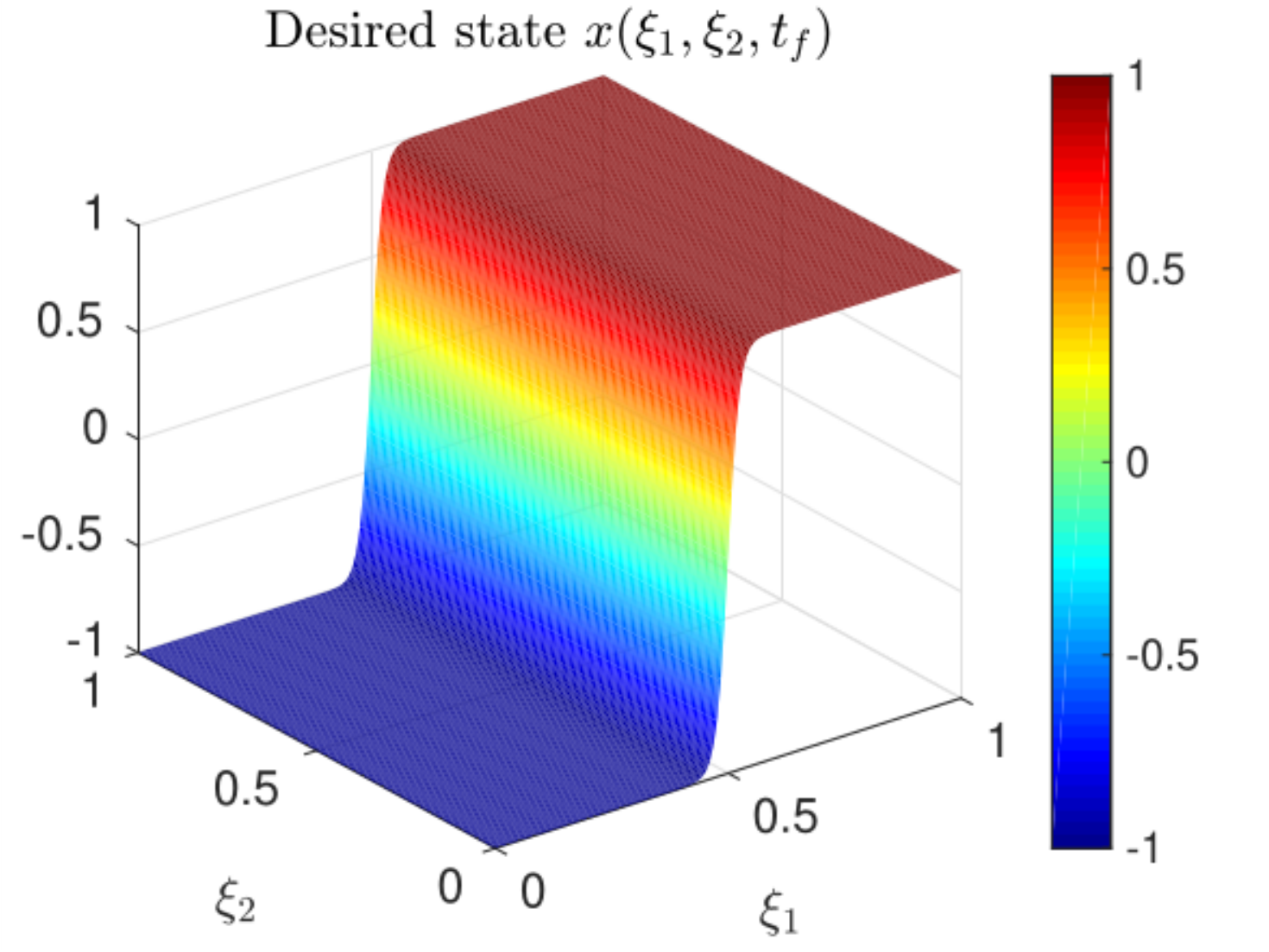}
		\includegraphics[scale=0.4]{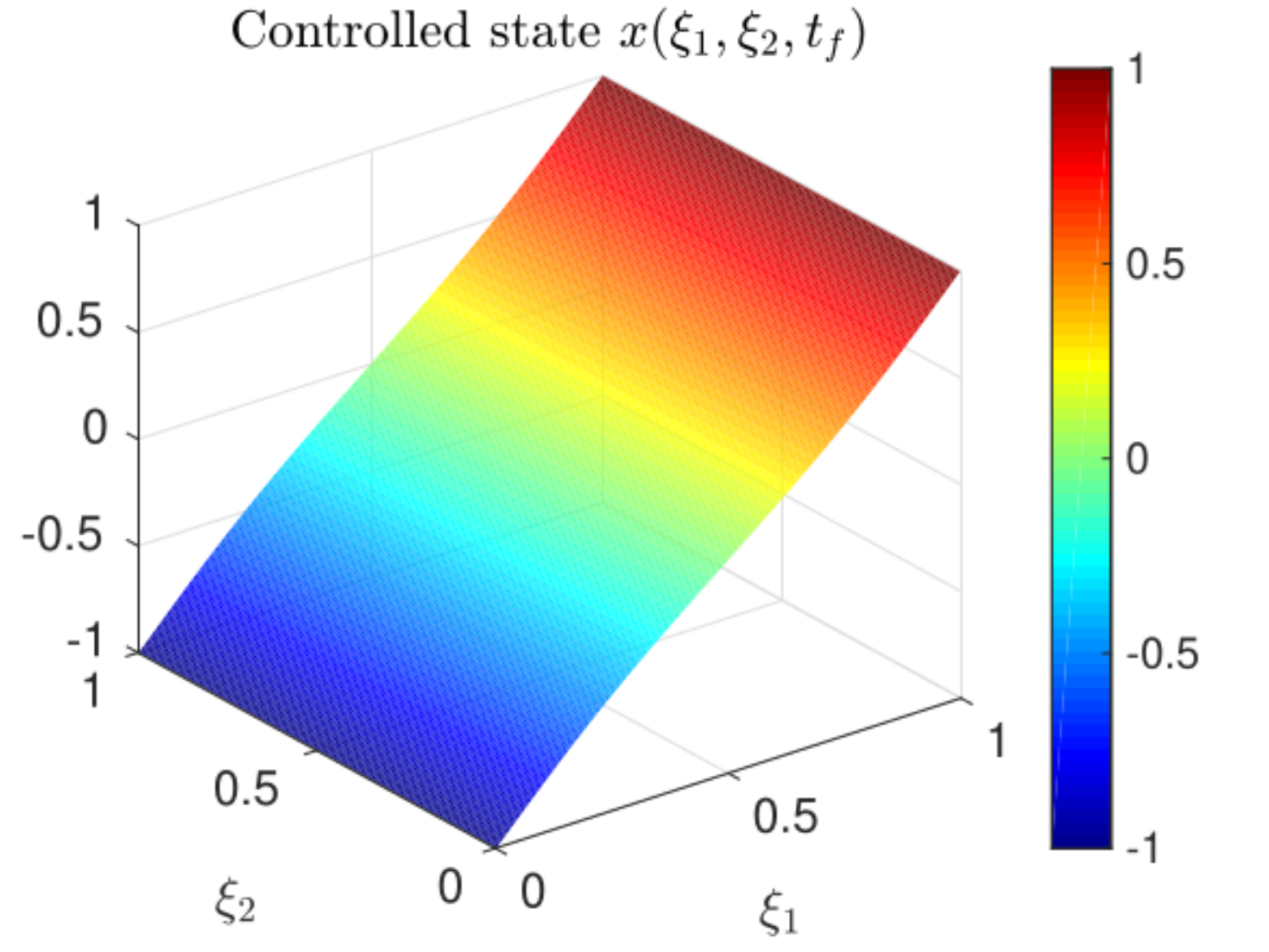}
		\caption{\emph{Left.} Desired state $x_d(\xi_1,\xi_2)$. \emph{Right.} Controlled state $x(\xi_1,\xi_2,t_f)$.}
		\label{fig:bound_des_vs_con_boundary}
	\end{center}
\end{figure}

In addition to the previous problem we also want to provide computational results for a heat equation equipped with a Neumann boundary control. In more detail, the PDE-constraint is given via 
\begin{align*} 
\frac{\partial}{\partial t} x(\xi,t) &= \Delta_{\xi} x(\xi,t)  && \text{in } \Omega \times (0,t_f), \\
\frac{\partial}{\partial \nu} x(\xi,t)&=u_1(t) && \text{on } \{0\}\times (0,1) \times (0,t_f), \\
\frac{\partial}{\partial \nu} x(\xi,t)&=u_2(t) && \text{on } \{1\}\times (0,1) \times (0,t_f), \\
\frac{\partial}{\partial \nu} x(\xi,t)&=u_3(t) && \text{on } (0,1)\times \{0\} \times (0,t_f), \\
\frac{\partial}{\partial \nu} x(\xi,t)&=u_4(t) && \text{on } (0,1)\times \{1\} \times (0,t_f), \\
x(\xi,0)&=x_0(\xi) && \text{in } \Omega,
\end{align*}
with domain $\Omega=(0,1)^2$ and boundary $\Gamma=\partial \Omega$. The observation $y(\cdot)$ is taken as in the previous case. However, the penalization parameter of the optimal control problem is set to $\beta=1$.
As for the case of an interior control, we consider tracking of a temporally constant state $x_d$ that is visualized in Figure \ref{fig:bound_des_vs_con_boundary}. The results shown in Figure \ref{fig:nk_iterates_ranks_boundary} indicate a similar behaviour to the the previously computed case. Namely, that eventually the quadratic convergence of the Newton method kicks in and that again relying on the nester approach saves many nonlinear iterations and thus more solves of the linearized space-time problem. Moreover, we do not observe an increase of the TT-ranks during the Newton-Kleinman iteration.

\subsection{Convection diffusion equation}
While the previous partial differential equation constraints are already quite challenging, we want to illustrate that our techniques are not limited to the case of the heat equation but that this approach is quite general. For this we here consider the convection diffusion equation given as
\begin{align*} 
\frac{\partial}{\partial t} x(\xi,t) &= \Delta_{\xi} x(\xi,t)+20 x_{\xi_1}(\xi,t)-10x_{\xi_2}(\xi,t)+200x(\xi,t)+ \sum_{i=1}^5 \chi_{\omega_i}(\xi) u_i(t)  
\end{align*}
defined on $ \Omega \times (0,t_f)$ and equipped with boundary condition $x(\xi,t)=0$  on  $\Gamma \times (0,t_f)$ and initial condition $x(\xi,0)=x_0(\xi)$ in $\Omega.$ We here provide a proof-of-concept implementation and use a simple finite difference scheme without additional upwinding. Note that in future research it is desirable to include a streamline upwind Petrov--Galerkin (SUPG) \cite{brooks1982streamline} or local projection stabilization \cite{ahmed2011discontinuous} scheme. We consider $t_f= 1$ and $\beta=10^{-4}$. Both control and observation domains are set in the case of the distributed heat equation. 

\begin{figure}[htb] 
	\begin{center}
		\includegraphics[scale=1]{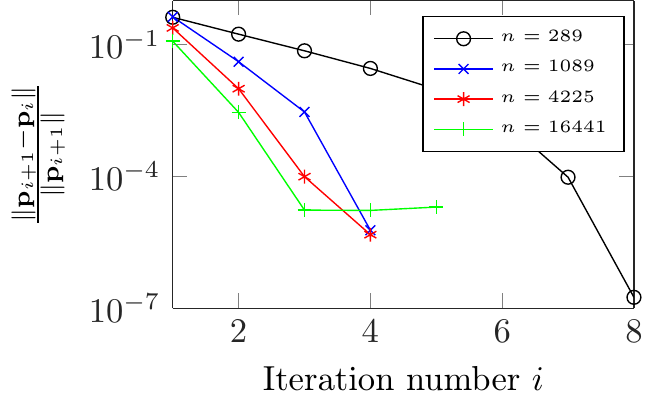}
		\includegraphics[scale=1]{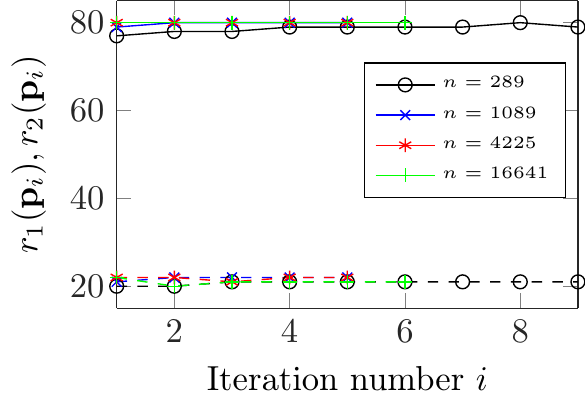}
		\caption{\emph{Left}. Change in Newton iterates. \emph{Right}. TT ranks of $\mathbf{p}_i$ during Newton iteration.}
		\label{fig:nk_iterates_ranks_interior_cd}
	\end{center}   
\end{figure}

\begin{figure}[htb]
	\begin{center}
		\includegraphics[scale=0.4]{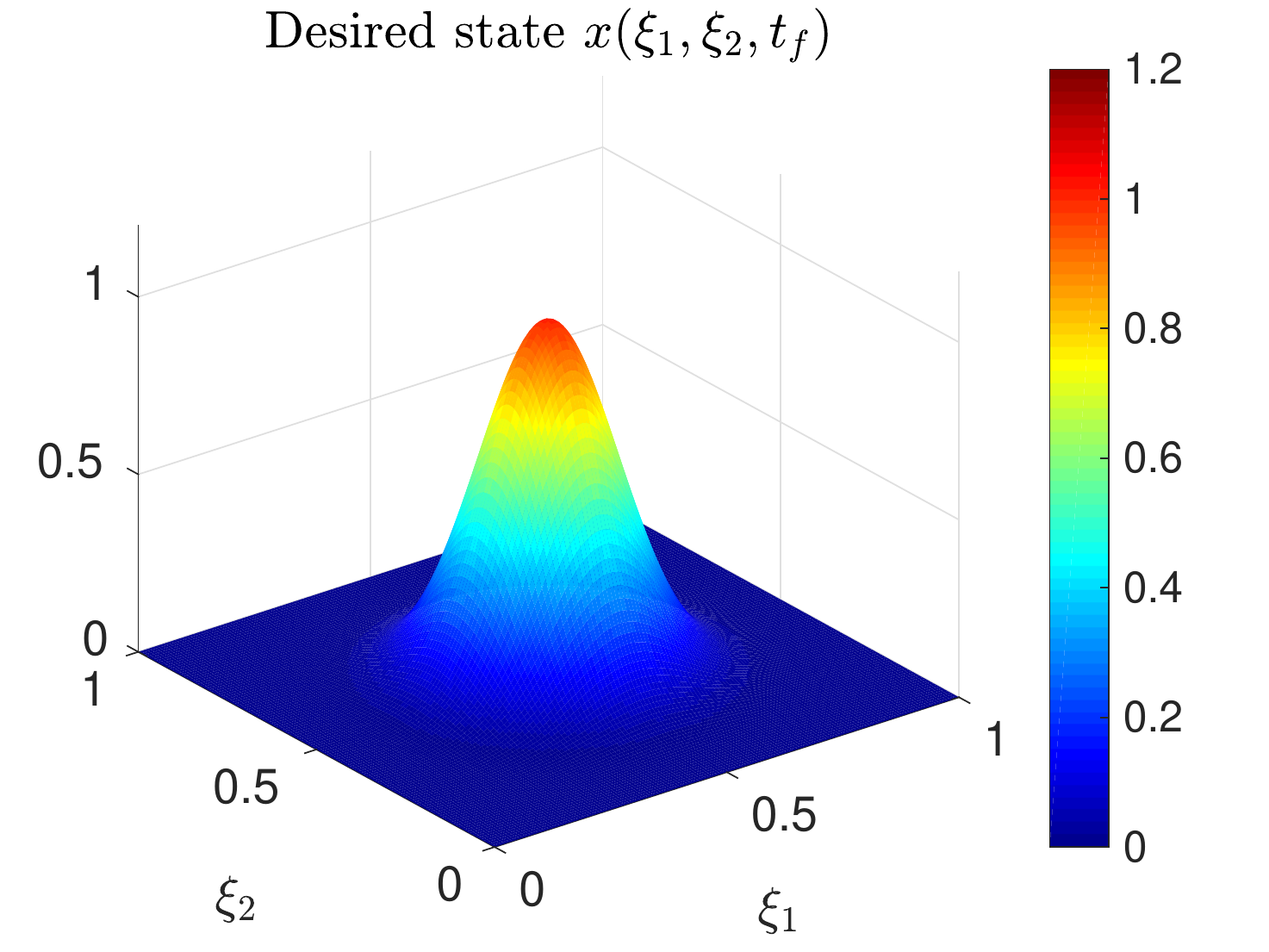}
		\includegraphics[scale=0.4]{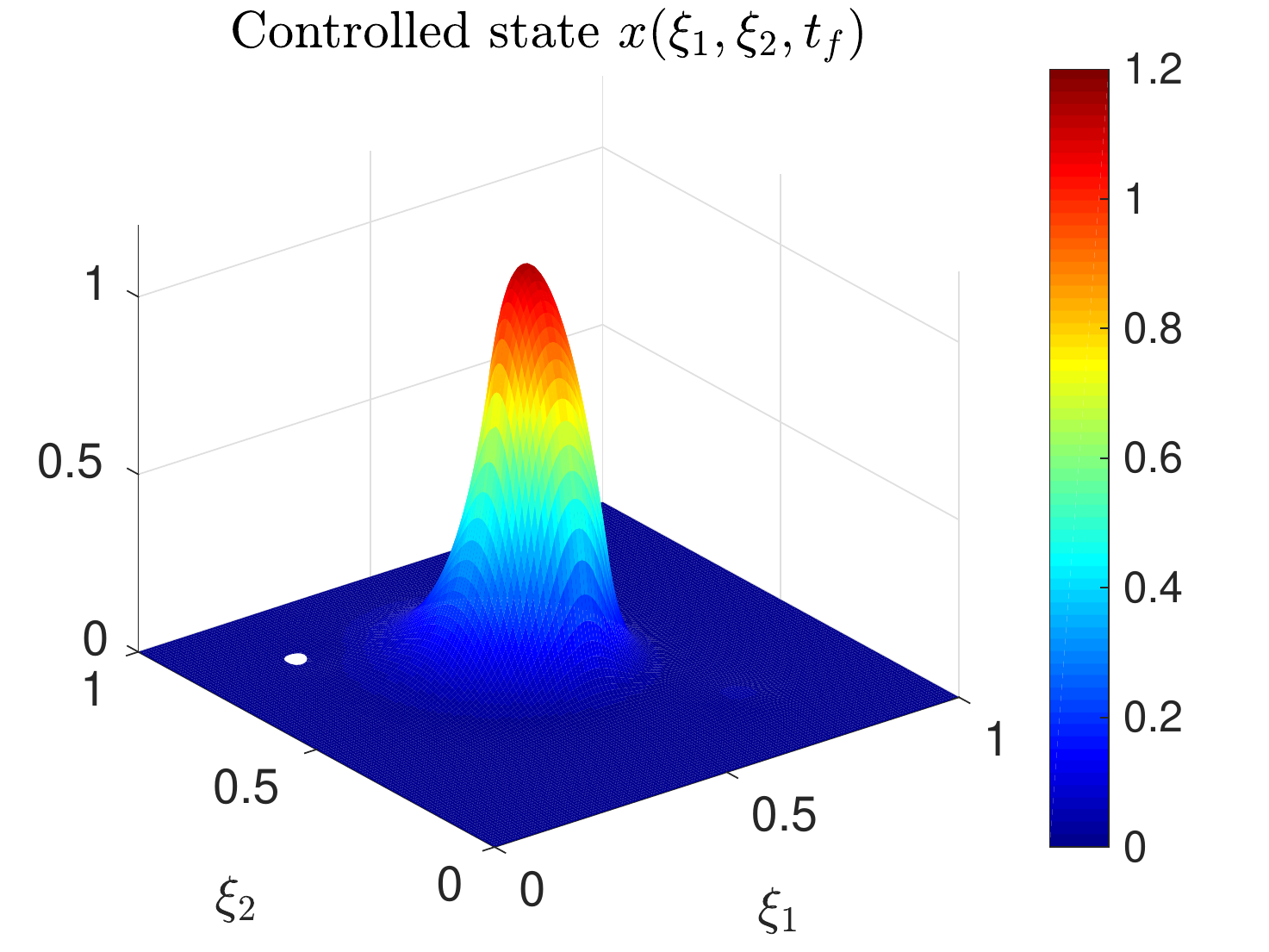}
		\caption{\emph{Left.} Desired state $x_d(\xi_1,\xi_2)$. \emph{Right.} Controlled state $x(\xi_1,\xi_2,t_f)$.}
		\label{fig:bound_des_vs_con_interior_cd}
	\end{center} 
\end{figure}

We show detailed results in Figure \ref{fig:nk_iterates_ranks_interior_cd} where again we employ a nested approach transferring the solution of the previous mesh level to the finer level as an initial guess. We then observe robust convergence but also note that the precision of $10^{-6}$ is not reached when the mesh contained $16641$ degrees of freedom. Additionally, we see that the ranks of the solution are similar to the ones needed for the heat equation control. In Figure \ref{fig:bound_des_vs_con_interior_cd} we see that the controlled state resembles the desired state rather well.
In order to compare qualitatively, we compared our solver to the MMESS\footnote{\url{https://www.mpi-magdeburg.mpg.de/projects/mess}} \cite{SaaKB19-mmess-2.0} library. As our approach is based on approximating the solution to the Riccati equation over the whole space-time domain we compared the storage for this tensor to the storage requirements within MMESS when the solution to all time-steps is stored in low-rank form. In this comparison we needed 230 MB of storage while the MMESS solutions required 21GB for a solution of average rank $157$ in each time-step.  

\section{Conclusion and outlook}
We here have provided a novel framework that has not been explored before, where we use an all-at-once linear low-rank solver for the differential Riccati equation. Namely, we posed the differential Riccati equation as a space-time problem with an outer Newton-Kleinman method handling the difficult nonlinearity. The inner linearized space-time problem is then highly structured and we designed a tailored low-rank tensor scheme that allows for the efficient solution to this problem. 

In the future we envisage to derive more tailored preconditioners for the subproblems within the AMEn method. This is especially crucial when more challenging PDEs are considered. Specifically for three-dimensional problem of this kind, it might be beneficial to use Tucker or Extended TT \cite{dk-qtt-tucker-2013} tensor formats. This would require a technical but straightforward algorithmic modification in the future.



\bibliographystyle{siam}
\bibliography{references}

\end{document}